\documentclass[12pt]{amsart} 
\usepackage{amsthm, amssymb,color,verbatim,mathpazo,graphicx}
\usepackage{fullpage,float} 
\usepackage[T1]{fontenc}
\newcommand{\lra}{\longrightarrow} 
\renewcommand{\le}{\leqslant} 
\renewcommand{\ge}{\geqslant} 
\renewcommand{\proof}{{\sc Proof. }}
\renewcommand{\P}{\mathbf P} 
\newcommand{\ZZ}{\mathcal Z}

\newcommand{\cA}{\mathcal A} 
\newcommand{\lcA}{\overline{\mathcal A}} 
\newcommand{\cB}{\mathcal B} 
\newcommand{\lcB}{\overline{\mathcal B}} 
\newcommand{\FFq}{{\mathbb F}_q} 
\newcommand{\bFFq}{{\overline{\mathbb F}}_q}

\newcommand{\ptc}{{\sc ptc}}
\newcommand{\ponc}{{\mathfrak P}}
\renewcommand{\phi}{\varphi} 

\newcommand{\Lower}{\mathbb L}
\newcommand{\Upper}{\mathbb U}
\newtheorem{Theorem}{Theorem}[section]
\newtheorem{Lemma}[Theorem]{Lemma}
\newtheorem{Proposition}[Theorem]{Proposition}
\newtheorem{Definition}[Theorem]{Definition} 
\newtheorem{Conjecture}[Theorem]{Conjecture} 
\begin{document} 
\title{On the Poncelet triangle condition over finite fields} 
\author{Jaydeep Chipalkatti} 

\maketitle

\bigskip 

\parbox{17cm}{ \small
{\sc Abstract:} Let $\P^2$ denote the projective plane over a finite
field $\FFq$. A pair of nonsingular conics $(\cA, \cB)$
is said to satisfy the Poncelet triangle condition if, considered as
conics in $\P^2(\bFFq)$, they intersect transverally and there exists a triangle 
inscribed in $\cA$ and circumscribed around $\cB$. It is shown in this
article that a randomly chosen pair of conics
satisfies the triangle condition with asymptotic probability $1/q$. We
also make a conjecture based upon computer experimentation which predicts this probability for
tetragons, pentagons and so on up to enneagons.} 

\bigskip \bigskip


\medskip 

\tableofcontents

\section{Introduction} 

\subsection{} \label{section.ponc.construction}
We begin by recalling Poncelet's closure theorem, which is one of
the most appealing results in classical projective geometry. 

Let $\P^2$ denote the projective plane over an 
algebraically closed field $\kappa$ of characteristic not $2$. 
Consider a pair of conics  $\cA$ and $\cB$ in $\P^2$ intersecting
transversally. Choose a point $P_1$ on $\cA$. Draw a tangent to $\cB$
from $P_1$, intersecting $\cA$ again at $P_2$. Now repeat the construction at $P_2$
to get a point $P_3$ on $\cA$, and then once again to get $P_4$. In general,
$P_4$ may not coincide with $P_1$; but if it does, then 
$P_1 P_2 P_3$ is a triangle inscribed in $\cA$ and circumscribed
around $\cB$. Such a triangle\footnote{In general there are two
  tangents to $\cB$ from $P_1$, from which we can opt for either one. If the
  triangle closes, then the other tangent automatically gets chosen at
  $P_3$.} will be called an $\cA \circ \cB$ triangle. 

\begin{figure}[H]
\includegraphics[width=8cm]{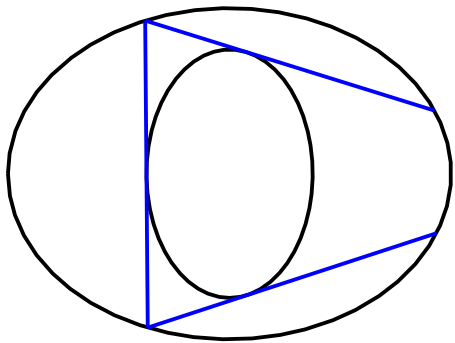} 
\caption{Conics failing the triangle condition} 
\label{diag:nptc} 
\end{figure}

\begin{figure}[H]
\includegraphics[width=8cm]{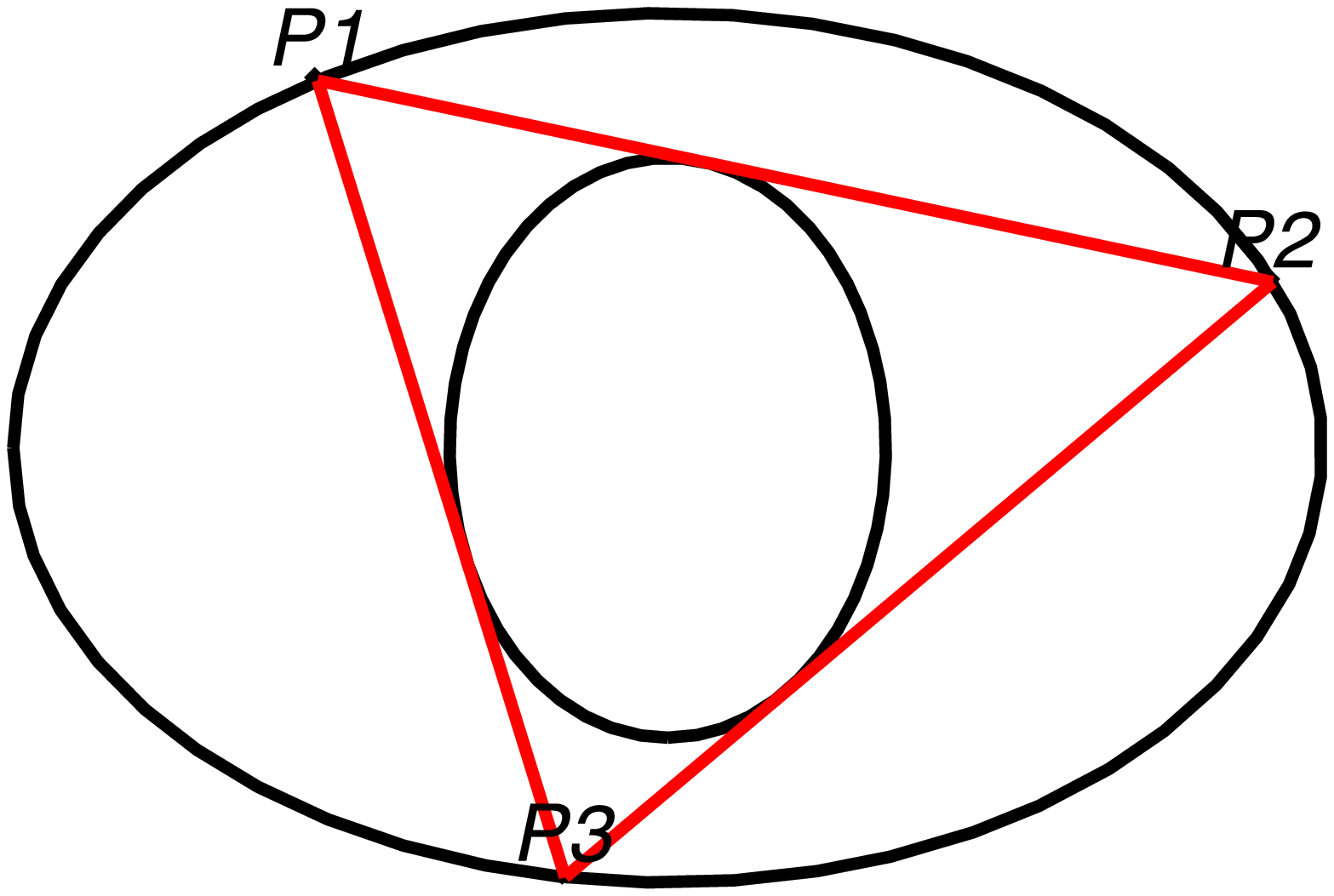} 
\includegraphics[width=8cm]{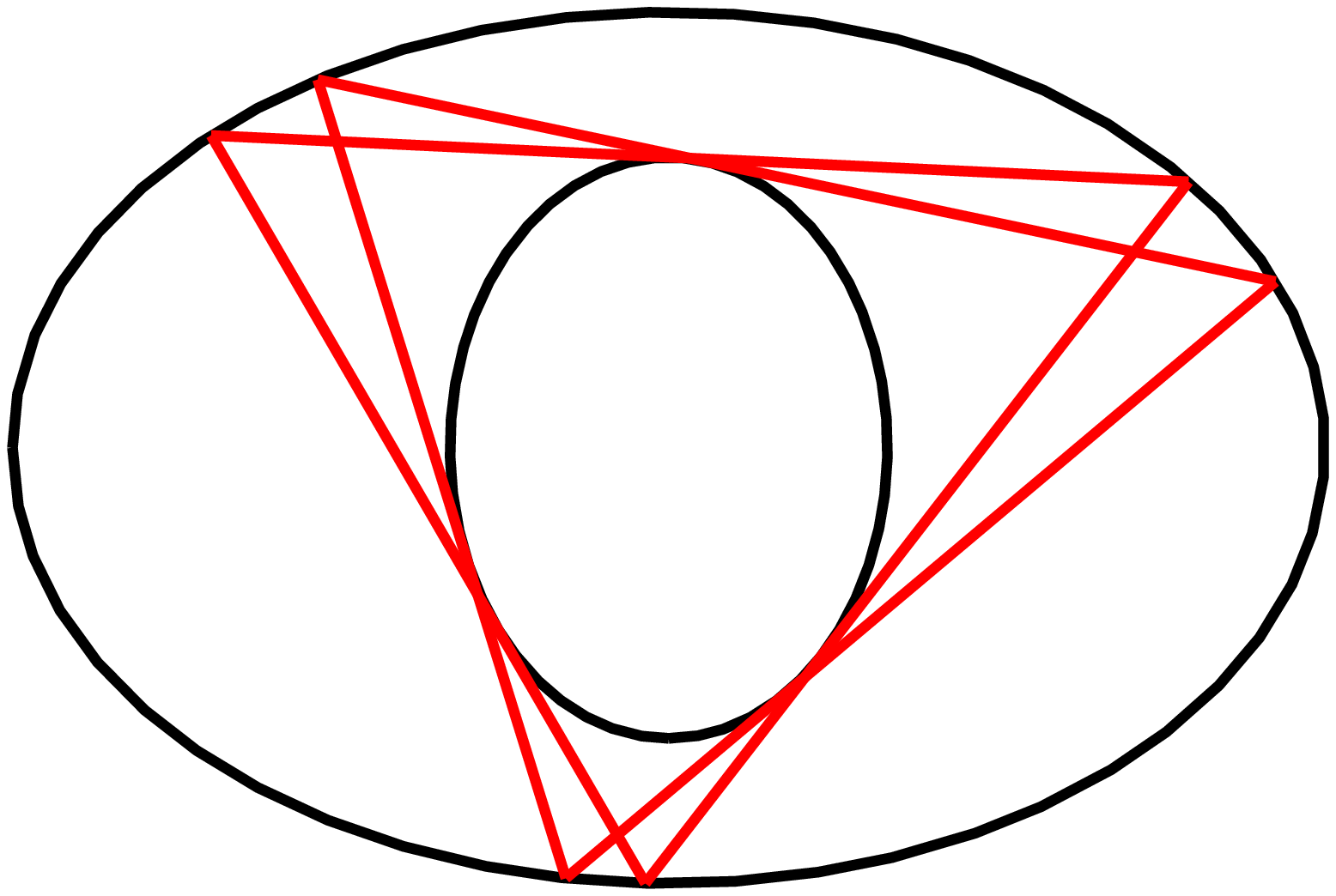}
\caption{Conics satisfying the triangle condition} 
\label{diag:ptc} 
\end{figure}

Now Poncelet's theorem says that if $P_4=P_1$ for \emph{some} 
choice of $P_1$, then the same is true of \emph{any} choice of
$P_1$. In other words, given $\cA$ and $\cB$, the problem of
constructing an $\cA \circ \cB$ triangle is poristic\footnote{According
  to John Playfair (writing in 1792): A Porism may be defined, 
A proposition affirming the possibility of finding such conditions as
will render a certain problem indeterminate, or capable of innumerable
solutions (source: Oxford English Dictionary).} in the sense that, 
it has either no solution or infinitely many
solutions (see Diagrams~\ref{diag:nptc} and \ref{diag:ptc}). The
former case is the norm and the latter the exception. There 
is no such triangle if the conics are~\emph{generally} situated; 
that is to say, they must be in geometrically special
position for the problem to be solvable. 

\subsection{} \label{section.defn.PTC} 
Now consider the plane $\P^2(\FFq)$ over a finite field
$\FFq$, where $q = p^r$ and $p \neq 2$. A conic $\cA \subseteq \P^2(\FFq)$
defines a conic $\lcA \subseteq \P^2(\bFFq)$ given by the same
equation. 

\begin{Definition} \rm \label{defn.ptc} 
We will say that a pair of nonsingular conics $(\cA, \cB)$ in $\P^2(\FFq)$ 
satisfies the Poncelet triangle condition (\ptc), if 
\begin{itemize} 
\item the conics $\lcA$ and $\lcB$ intersect transversally (i.e., in four distinct points), 
\item 
There exists an $\lcA \circ \lcB$ triangle. 
\end{itemize} 
\end{Definition}

\medskip 

Since (\ptc) is a nontrivial geometric condition on the pair,
it is natural to ask how frequently one can expect it to hold. 
The main result of this paper (Theorem~\ref{main.theorem}) can be
paraphrased as saying that, 

\smallskip 

\centerline{The proportion of conic pairs satisfying (\ptc) is
asymptotically $\frac{1}{q}$.} 

\medskip 

In other words, the probability that a randomly chosen conic pair
satisfies (\ptc) is approximately $\frac{1}{q}$. The 
actual statement of the theorem gives an upper and a lower bound for
this proportion. 

\medskip 

Two clarifications are in order: 
\begin{itemize} 
\item 
The conics have $4$ common points in $\P^2(\bFFq)$, and either $0, 1,
2$ or all $4$ of them will be in $\P^2(\FFq)$. 
\item 
The definition of (\ptc) by itself does not require that there be an $\cA \circ
\cB$ triangle. However, there do exist such triangles
when $q$ is sufficiently large (see section~\ref{section.ABtriangle}). 
\end{itemize}

\subsection{} Poncelet's theorem overlaps several
areas of mathematics, and as such the literature associated to it is
very large. The article by Bos et.~al.~\cite{Bos} is a masterly survey of the historical
development of the theorem. It contains an account of
Poncelet's own proof, as well as Jacobi's proof using elliptic
functions. Halbeisen and Hungerb{\"u}hler~\cite{HH} give another proof
using Pascal's theorem. One can also find a wealth of material in the 
treatises by Dragovi{\'c}-Radnovi{\'c}~\cite{DR} and
Flatto~\cite{Flatto}. The preprint by 
Hungerb{\"u}hler and Kusejko~\cite{HK} contains an  
interesting discussion of Poncelet's theorem for projective planes
over prime fields. We refer the reader to Coxeter~\cite{Coxeter}
and~Hirschfeld~\cite[Ch.~7]{Hirschfeld} for standard facts about conics
in projective planes. 

\medskip 

{\sc Acknowledgement:} I am thankful to Sudhir Ghorpade and Keith
Mellinger for some helpful correspondence. 

\subsection{} Although a complete proof of Poncelet's theorem will 
not be reproduced here, we enclose a summary of the now-classic
Griffiths-Harris proof~\cite{GH} for the reader's interest. 
Assume the base field to be algebraically
closed of char $\neq 2$, and that $\cA, \cB$ intersect transversally. Let $\cB^* \subseteq
(\P^2)^*$ denote the dual conic consisting of tangent lines to
$\cB$. Consider the subvariety
$ E \subseteq \cA \times \cB^*$ given by 
\[ E = \{ (P,T): \text{$T$ is a tangent to $\cB$ passing through $P$} \}. \] 
The projection morphism $E \lra \cA$ is a double cover branched over the four
points in $\cA \cap \cB$. It follows by the Riemann-Hurwitz formula
that $E$ is an elliptic curve. The function 
$(P_i, \overline{P_i P_{i+1}}) \lra (P_{i+1}, \overline{P_{i+1}
  P_{i+2}})$, from section \ref{section.ponc.construction} corresponds to a
translation 
\[ E \lra E, \quad z \lra z + \tau, \] 
by some constant $\tau \in E$. Now $P_4=P_1$, iff $\tau$ is a $3$-torsion point of
$E$. But $\tau$ depends only on the relative positions of $\cA$
and $\cB$, and hence $P_4 = P_1$ is true either for no $P_1$ or for
all $P_1$. \qed 

\smallskip 

The argument remains unchanged if $3$ is replaced by any
$n$. Thus, if there exists an $n$-gon inscribed in $\cA$ and
circumscribed around $\cB$, then there exists one starting from any
point in $\cA$. Although the main result of this paper applies only to 
triangles, we propose a conjecture about the next few values of $n$ 
(see section~\ref{section.higher_n}). 

\subsection{} 
Even if the pair $(\cA, \cB)$ satisfies~(\ptc), it may happen over a
finite field that no tangent can be drawn to $\cB$ from some choices 
of $P_1$ in $\cA$. (This will be the
case if the polar line of $P_1$ with respect to $\cB$ does not
intersect $\cB$ in an $\FFq$-rational point.) However, if such a 
tangent does exist, then one can complete an $\cA \circ \cB$
triangle. Examples of either phenomenon will be given in section 
\ref{example.ptc}.

\section{The Main Theorem} 

Assume that $\text{char}(q) \neq 2, 3$. 
Let $\Psi$ denote the set of conic pairs $(\cA, \cB)$ in $\P^2(\FFq)$,
such that $\lcA, \lcB$ intersect transversally. Let $\Gamma$ denote the subset of pairs
satisfying (\ptc). 

\begin{Theorem} \rm \label{main.theorem} 
With notation as above, 
\[ \frac{q-16}{q \, (q+1)} \le \frac{|\Gamma|}{|\Psi|}
\le \frac{q+5}{(q-2)(q-3)}. \] 
\end{Theorem} 

One can think of a conic pair in $\Psi$ as being a candidate 
for satisfying (\ptc). According to the theorem, the probability that it actually does so 
is $\frac{1}{q} + O\left(\frac{1}{q^2}\right)$.

\subsection{} 
Our main tool will be an algebraic criterion due to Cayley
for (\ptc) to hold (see~\cite{GH}). 
Let $[x,y,z]$ be homogeneous coordinates in $\P^2$. Let $\cA, \cB$, 
which are assumed to intersect transversally in $\P^2(\bFFq)$, 
respectively correspond to symmetric $3 \times 3$ matrices $A, B$. 
Write $\Delta = \det (t \, A + B)$, where $t$ is an indeterminate. Now
consider a formal Maclaurin series expansion 
\[ \sqrt{\Delta} = H_0 + H_1 \, t + H_2 \, t^2  + \dots,  \] 
where the $H_i$ are functions of entries in $A$ and $B$. Then we have
the following criterion: 

\begin{Proposition}[Cayley] \rm The pair $(\cA, \cB)$ satisfies (\ptc), if and only
  if $H_2=0$. 
\end{Proposition} 

In the context of the proof-sketch above, $E$ has an affine
model given by the equation $u^2 = \Delta$ in the variables $t,u$. Now
the criterion is proved by an explicit calculation which detects the inflection
points of $E$ (see~[loc.~cit.]).

\subsection{A sample calculation} \label{section.sample.calc} 
We will begin by determining the (\ptc)-pairs in a special
case. Most of the ideas needed for the main theorem are already present in
this calculation. 
For simplicity, assume that $q$ is a prime number $\ge 7$. 
Consider the pencil of conics in $\P^2(\FFq)$ passing through the points 
$[1,0,0], [0,1,0], [0,0,1], [1,1,1]$. Each nonsingular conic in the
pencil may be written as 
\[ C_\alpha: \alpha \, xy + (1-\alpha) x \, z - yz =0, \] 
for some $\alpha \in \FFq \setminus \{0,1\}$. It corresponds to the
symmetric matrix 
$\left[ \begin{array}{rrr} 0 & \alpha & 1-\alpha \\ 
\alpha & 0 & -1 \\ 1-\alpha & -1 & 0 \end{array} \right]$. Let 
$\cA = C_r$ and $\cB = C_s$ for some $r,s \neq 0, 1$; the number of
such pairs is $(q-2) \, (q-3)$. Now consider the subset 
\[ \ponc = \{(r,s): H_2(r,s)=0, \; \text{and} \; r, s \neq 0,1 \}, \]  
of conics satisfying (\ptc). We will determine the size of $\ponc$. 

A straightforward calculation shows that $\Delta = (r \, t+s) \, (r \,
t+s-t-1) \, (t+1)$, 
and\footnote{The denominator of $H_2$ is $[s (s-1)]^{\frac{3}{2}}$, which can be ignored.} 
\begin{equation} 
H_2(r,s) = r^2 + (6 \, s^2 -4 \, s^3 - 4 \, s) \, r + s^4. 
\label{ponc3} \end{equation} 
Considered as a quadratic in $r$, its discriminant is 
\begin{equation} \delta = (6 \, s^2 -4 \, s^3 - 4 \, s)^2 - 4 \, s^4 = 
\underbrace{16 \, s^2 \, (s-1)^2}_{e(s)} \, \times 
\underbrace{(s^2-s+1)}_{f(s)}. 
\label{eqn.delta} \end{equation}  
Thus, for a given value of $s$, the equation $H_2(r,s)=0$ has one root
in $r$ if $f(s)=0$, two distinct roots if $f(s)$ is a nonzero 
square in $\FFq$, and no roots otherwise. 

\medskip 

\noindent {\bf Claim-1:} The set $\ponc$ has cardinality $q-5$. 

\smallskip  

\noindent Given claim-1, the proportion of (\ptc)-pairs in the pencil
is 
\[ \frac{q-5}{(q-2) (q-3)} \simeq \frac{1}{q}.  \] 

\smallskip 

To prove claim-1, consider the set 
$S = \{s \in \FFq \setminus \{0,1\}: \text{$f(s)$ is a nonzero
square}\}$. Let $\left(\frac{m}{q} \right)$ denote the Legendre symbol. 

\medskip 

\noindent {\bf Claim-2:} The size of $S$ is given by 
\[ |S| = \begin{cases} 
\frac{q-7}{2} & \text{if $\left(\frac{-3}{q}\right) = +1$}, \\ 
\frac{q-5}{2} & \text{if $\left(\frac{-3}{q}\right) =-1$}. 
\end{cases} \]

Claim-1 follows from claim-2. Indeed, suppose that
$\left(\frac{-3}{q}\right) = +1$. Then $S$ contributes $2 \times
\frac{q-7}{2} = q-7$ pairs to $\ponc$. Furthermore, 
$f(s) = (s-1/2)^2 + 3/4 =0$ has two roots in $\FFq$, each of which
contributes one pair, making up the total of $q-5$. The
other case follows by a similar argument. 

\medskip 

To prove claim-2, write $(s-1/2)^2+3/4 = y^2$ for $s \in S$. From 
\[ \underbrace{(y-s+1/2)}_a \, \underbrace{(y+s-1/2)}_{3/4a}  =
3/4, \] 
we get $s = \frac{3+4a-4a^2}{8a}$ for some $a \neq 0$. Since $s \neq
0,1$, we must have $a \neq \pm 1/2, \pm 3/2$. Altogether this excludes
five values of $a$. If $\frac{1}{2} (1 \pm \sqrt{-3}) \in \FFq$, then
two more values (namely $a=\pm \frac{\sqrt{-3}}{2}$) are excluded,
since $f(s)=0$ is disallowed. Since $a$ and $-3/4a$ lead to
the same $s$-value, we must divide the number of possible $a$-values by $2$. Thus we get 
either $\frac{q-7}{2}$ or $\frac{q-5}{2}$, according to whether $\sqrt{-3}$
does or does not belong to\footnote{This condition can be made explicit using 
the quadratic reciprocity theorem (see~\cite[Ch.~5]{IR}). 
We have $\left( \frac{-3}{q} \right) = +1$ (resp.~$-1$) if $q \equiv
1, 7$ (resp.~$5,11$) mod $12$.} $\FFq$. \qed 

\subsection{Example} \label{example.ptc} 
Let $q=43, \cA = C_{11}$ and $\cB = C_{36}$. It is easy to verify that
$H_2=0$, so that (\ptc) holds. Let $P_1=[1,17,34]$ on $\cA$. The 
polar line to $\cB$ with respect to $P_1$ is $x+18y+5z=0$. It 
intersects $\cB$ in the two points $R_1 = [1,32,5]$ and $R_1' = [1,40,2]$. If we choose 
$R_1$, then $P_1R_1$ is a tangent to $\cB$ which intersects $\cA$ again at $P_2 =
[1,36,3]$. Repeating the construction at $P_2$ leads to $P_3 =
[1,24,28]$, and then back to $P_1$. This gives an $\cA \circ \cB$
triangle.  By contrast, if $P_1=[1,9,12]$, its polar line with respect 
to $\cB$ is $x+32y+13z=0$, which does not intersect $\cB$. The 
construction cannot proceed any further, and there 
is no $\cA \circ \cB$ triangle with $[1, 9, 12]$ as a vertex. (Of course,
there will exist an $\lcA \circ \lcB$ triangle.) 

A choice of $P_1$ on $\cA \cap \cB$ will lead to a 
degenerate triangle, with two coincident vertices. For instance, with 
the same conics as above, let $P_1 = [0,1,0]$. 
The tangent to $\cB$ through $P_1$ intersects $\cA$ again 
at $P_2=[1,20,36]$. Now the tangent to $\cB$ through $P_2$ is 
$x+14y+34z=0$, which is also the tangent to $\cA$ at 
$P_2$. Hence $P_3$ coincides with $P_2$. 

\subsection{} \label{section.invariance} 
It is easy to check directly that
\[ H_2(1-r,1-s) = H_2(r,s), \quad e(1-s) =e(s), \quad f(1-s)=f(s). \] 
In geometric language, interchanging $y$ and $z$ defines an involution on
the pencil which interchanges $C_\alpha$ and $C_{1-\alpha}$. This ${\mathbf Z}_2$-invariance
will play a small role later. 
\section{The Dickson classification} 

In this section we will complete the proof of the main theorem. In
outline, the strategy is to decompose $\Psi$ into a union of pencils, and
estimate the proportion of (\ptc)-pairs in each pencil. 

\subsection{} 
Two quadratic forms $F(x,y,z), G(x,y,z)$ form a pencil 
\[ \pi = \{\eta \, F + G=0: \eta \in \P^1(\FFq) \} \subseteq \P^2(\FFq). \] 
By convention, $\eta = \infty$ corresponds to $F=0$. All such
pencils have been classified up to projective automorphisms by
Dickson~\cite{Dickson}; this classification is also described in a
table on page 175 of Hirschfeld~\cite{Hirschfeld}. There are
altogether $20$ isomorphism 
classes, but $15$ of them are fortunately disqualified for at least one of the
following reasons: 
\begin{itemize} 
\item 
Every member of the pencil is singular (e.g., $(1)$-st entry in
the table). 
\item 
The generators of the pencil do not intersect transversally in
$\P^2(\bFFq)$ (e.g., $(4)$-th entry). 
\item 
The pencil can only occur in characteristic $2$ (e.g., $(7)$-th
entry). 

\end{itemize} 
This leaves five isomorphism classes, to be called class $(i)$ for $i=3,14, 16,
18, 19$, corresponding to their positions in the table. 

For class $(3)$, the generators are 
$F = xy, G = z^2+yz+xz$. They intersect in four 
$\FFq$-rational points, namely $[0,1,0], [0,1,-1], [1,0,0], [1,0,-1]$. The pencil 
in section~\ref{section.sample.calc} is of this class. 

For class $(14)$, the generators are 
\[ F = xy, \quad G = y^2+yz + xz + e \, z^2, 
\] 
where $e \in \FFq$ is any element such that the polynomial $T^2+T+e$ is  
irreducible over $\FFq$. 

For class $(16)$, the generators are 
$F = xy, G = e_1 \, x^2 + e_2 \, y^2 + xz + yz + z^2$, 
where $e_1, e_2$ are such that $T^2 + T+e_1, T^2 + T+e_2$ are irreducible. 

For class $(18)$, the generators are 
$F = y^2 - xz, G = x^2 + b \, y^2 + c \, xy + yz$, 
where $b,c$ are such that 
\begin{equation} g(T) = T^3 + bT^2 + cT +1, \label{gT} \end{equation} 
is irreducible. 

For class $(19)$, the generators are 
$F = x^2 - \nu\, y^2, G = z^2 - \rho \, y^2 + 2 \, \sigma \, xy$, 
where $\nu$ and $\rho^2-4 \, \nu \, \sigma^2$ are non-squares. 

Notice that $F=0$ is a singular conic for all classes except
$(18)$. Each of the other four classes has three singular members
corresponding to $\eta_1, \eta_2, \eta_3$, where
$\eta_1, \eta_2 \in \bFFq$ are the roots of the quadratic equation 
$\text{discrim}(\eta \, F + G)=0$, and $\eta_3= \infty$. For class
$(18)$, this equation is a cubic in $\eta$, which is in fact identical
to $g(\eta)=0$. Because of this small anomaly, one will have to make a separate
argument for class $(18)$. 
\subsection{} 
We will say that a pencil $\pi$ is~\emph{eligible} if it belongs to
any of these five isomorphism classes. Let $\Psi_\pi$ denote the set of 
pairs of nonsingular conics in $\pi$, and let
$\Gamma_\pi$ denote the subset of pairs satisfying (\ptc). Write 
\[ \Lower = \frac{q-16}{q \, (q+1)}, \quad 
\Upper = \frac{q+5}{(q-2) \, (q-3)}. \] 
for the lower and upper bounds in the main theorem. 

\begin{Proposition} \rm \label{prop.estimate.pi} 
For every eligible pencil $\pi$, we have 
\[ \Lower \le 
\frac{|\Gamma_\pi|}{|\Psi_\pi|} \le \Upper. \] 
\end{Proposition} 
Given the proposition, the main theorem follows immediately. 
We have decompositions 
$\Psi = \bigcup\limits_\pi \Psi_\pi$ and $\Gamma = \bigcup\limits_\pi
\Gamma_\pi$, quantified over eligible pencils. Then 
\[ \Lower \le \frac{|\Gamma|}{|\Psi|} = 
\frac{\sum_\pi |\Gamma_\pi|}{\sum_\pi
  |\Psi_\pi|} \le \Upper. \] 
Here we have used the elementary inequality 
\[ \min\limits_i \left\{\frac{a_i}{b_i} \right\} \le 
\frac{a_1 + a_2 + \dots}{b_1 + b_2 + \dots} \le 
\max\limits_i \left\{\frac{a_i}{b_i} \right\}. \] \qed  

\subsection{} 
The central idea behind the proposition is that the
structure of $H_2$ and $\delta$ for any eligible pencil is similar to the one
in the sample calculation, which allows us to make a qualitative estimate
along the lines of claim-2. We will break down the argument in a
couple of lemmas. Let $\cA, \cB$ respectively
correspond to $A = r \, F+G, B = s \, F + G$. 
\begin{Lemma} \rm Assume that $\pi$ is
  of any class except $(18)$. Then 
\begin{itemize} 
\item
The polynomial $H_{2,\pi}(r,s)$ is of degree $2$ in $r$, and degree $3$
in $s$. 
\item 
Its $r$-discriminant $\delta_{\pi}$ is of the form 
$e_\pi(s) \times f_\pi(s)$, where $e_\pi(s)$ is the square of a
quadratic polynomial, and $f_\pi(s)$ is a quadratic polynomial
which is not the square of a linear polynomial. 
\end{itemize} 
\label{prop.delta} \end{Lemma} 
For instance, for the class $(14)$ pencil, 
\[ \delta_\pi = \underbrace{16 \, (e \, s^2 - s +1)^2}_{e_\pi(s)} \times 
\underbrace{(e^2 \, s^2 - e \, s - 3 \, e +1)}_{f_\pi(s)}. \] 
If $f_\pi(s)$ were to be the square of a linear
polynomial, its $s$-discriminant 
\[ e^2 - 4 \, e^2 \, (1-3 \, e) = 3 \, e^2 \, (4 \, e-1) =0, \] 
implying that $e=0$ or $\frac{1}{4}$. But then $T^2 + T+e$ cannot be
irreducible in either case. 

\bigskip 

\proof Let $\overline{\pi}$ denote the corresponding pencil (defined by the same
generators) in $\P^2(\bFFq)$. Its base locus consists of a quadruple of
non-collinear points, and any two such quadruples can be taken to each
other via an automorphism of $\P^2(\bFFq)$. Thus all such pencils are
isomorphic over $\bFFq$, and we can obtain $H_{2,\pi}$ and $\delta_\pi$
by transforming the corresponding expressions ~(\ref{ponc3}), (\ref{eqn.delta})
from section~\ref{section.sample.calc}. Since the singular members must correspond, 
the two pencils are related by the affine substitution
$\alpha = \frac{\eta-\eta_1}{\eta_2-\eta_1}$, so that 
$\eta = \eta_1, \eta_2, \infty$ respectively map to $\alpha=0,1,\infty$. 
But then the same substitution on $r,s$ transforms $H_2, e(s)$ and
$f(s)$ into $H_{2,\pi}, e_\pi$ and $f_\pi$. Since all degrees are
preserved, and the property of being a square or
a non-square is likewise preserved, we have the result. \qed 

\bigskip 

The coefficients of the substitution are in $\bFFq$, and not
necessarily in $\FFq$. However, notice that 
$1- \alpha = \frac{\eta - \eta_2}{\eta_1 - \eta_2}$, which is
$\alpha$ with $\eta_1, \eta_2$ interchanged. Now the invariance
in section~\ref{section.invariance} 
implies that the coefficients of $H_{2,\pi}, e_\pi, f_\pi$ are
symmetric in $\eta_1, \eta_2$, and hence lie in $\FFq$. 
Such an argument will not work on class $(18)$, since one would need 
a fractional linear transformation to move $\alpha=\infty$ to a finite point $\eta_i$. 
\subsection{} 
We can now estimate the size of
$\Gamma_\pi$. The idea, as before, is to consider how often $f_\pi(s)$
is a square. Let $\phi(s) \in \FFq[s]$ be any 
quadratic polynomial which is not the square of a linear polynomial, and let
\[ \ZZ_\phi= \{s \in \FFq: \text{$\phi(s)$ is a square in $\FFq$} \}. \] 

\begin{Lemma} \rm 
With notation as above, 
\[ \frac{q-1}{2} \le |\ZZ_\phi| \le \frac{q+5}{2}. \] 
\label{lemma.Zf.estimate} 
\end{Lemma} 
\proof We will be brief, since the argument is similar to claim-2. 
(But the values $s=0,1$ are no longer
disallowed.) Say $\phi(s) = u_0 \, s^2 + u_1 \, s + u_2$. If $u_0$ is a
square, then dividing by it leaves $\ZZ$ unchanged, hence we may
write $f(s) = (s-b)^2 + c$ for some $c \neq0$. Then, as in claim-2, each element in
$\ZZ_\phi$ is of the form $s=\frac{c+2 \, b \, a -a^2}{2 \, a}$ for some $a
\neq 0$. Now $a, -c/a$ lead to the same $s$-value, and hence we get 
$|\ZZ_\phi| = \frac{q+1}{2}$ or $\frac{q-1}{2}$ depending on whether
$\sqrt{-c}$ does or does not belong to $\FFq$. In any case, 
\[ \frac{q-1}{2} \le |\ZZ_\phi| \le \frac{q+1}{2}. \] 

If $u_0$ is not a square, then consider $\tilde{\phi}(s) = \phi(s)/u_0$. Then $\phi(s)$
is a square iff $\tilde{\phi}(s)$ is a non-square, unless they are both zero. Applying the
earlier estimate to $\tilde{\phi}$ and taking complements, we get 
\[ \frac{q-1}{2} \le |\ZZ_\phi| \le \frac{q+5}{2}. \]
This proves the lemma. \qed

\subsection{} \label{section.Gammaf.estimate} 
Now let $\pi$ be an eligible pencil, not of class $(18)$, with singular
members $\eta_1, \eta_2$. Depending on its structure, 
it may happen that both $\eta_i$ belong to $\FFq$ or neither of
them does. 

Since an element in $\ZZ_{f_\pi}$ can
contribute at most two pairs to $\Gamma_\pi$, we have 
$|\Gamma_\pi| \le 2 \, |\ZZ_{f_\pi}|$.  It remains to find a
lower bound. We get only one $r$-value if $f_\pi(s)=0$. Since there are at
most two roots of $f_\pi(s)$ in $\FFq$,
this means a loss of at most two pairs. Moreover, at most $2 \times
2=4$ pairs may be lost because either $r$ or $s$ equals $\eta_i$. Thus 
$|\Gamma_\pi| \ge 2 \, |\ZZ_{f_\pi}| - 6$. 
Combining with the previous lemma, 
\begin{equation} 
 q-7 \le |\Gamma_\pi|\le q+5,  
\label{estimate1} \end{equation} 
for all eligible pencils except those of class $(18)$. 

\subsection{} The argument for class $(18)$ is a little more
intricate, but not different in substance. Calculating directly from
the pencil generators, we get 
$H_{2,\pi}(r,s) = h_0 \, r^2 + h_1 \, r + h_2$, where 
$h_0 = 3s^4 + 4 \, b \, s^3 + 6 \, c \, s^2 + 12 \, s +4 \, b-c^2$, 
and $h_1, h_2$ are polynomials in $s$ which need not be 
written down explicitly. Thus $h_0 $ is nonzero for all but at most $4$ values of 
$s$. The $r$-discriminant of $H_{2,\pi}$ is 
\[ \delta_\pi = h_1^2 - 4 \, h_0 \, h_2  = 
\underbrace{16 \, (s^3 + b \, s^2 + c \, s +1)^2}_{e_\pi(s)}
\times 
\underbrace{[(b^2-3 \, c) \, s^2 + (bc-9) \, s +
  (c^2-3b)]}_{f_\pi(s)}. \] 
In minor contrast to the earlier cases, $e_\pi$ is the square of a
cubic, and $f_\pi$ is of degree \emph{at most} $2$. 

We claim that the coefficients of $s^2$ and $s$ in $f_\pi$ cannot
vanish simultaneously. If they did, then $b^2 = 3c, bc =9$ and $b,c
\neq 0$ would together imply that $b =
3 \, \omega, c = 3 \, \omega^{-1}$, where $\omega$ is a cube-root of unity. But 
then the polynomial $g(T)$ from~(\ref{gT}) has $-\omega$ as a root, which 
contradicts its irreducibility. 
Moreover, $f_\pi(s)$ cannot be the square of a linear form. If it were, then its 
$s$-discriminant $-3b^2c^2+12b^3+12c^3-54bc+81=0$. But then $g(T)$
cannot be irreducible, since its $T$-discriminant is 
$b^2c^2-4b^3-4c^3+18bc-27=0$. (This follows from Dickson's
irreducibility criterion for cubics over finite fields --
see~\cite[Theorem~3]{Dickson2}.) 

\subsection{} 
If $f_\pi$ is a quadratic (i.e., if $b^2 \neq 3c$), then apply lemma 
\ref{lemma.Zf.estimate}. The argument in
section~\ref{section.Gammaf.estimate} also goes through, except that
we may lose at most $3 \times 3 =9$ pairs due to singular values. If
$h_0(s)=0$, then we get only one $r$-value, hence we may lose $4$ more
pairs this way. Thus 
\begin{equation} q-16 \le |\Gamma_\pi| \le q+5.  
\label{estimate2} \end{equation} 
If $f_\pi(s)$ is linear (i.e., if $b^2=3c$), then it is a square 
for $\frac{q+1}{2}$ values of $s$, which gives 
\begin{equation}  q-14 \le |\Gamma_\pi| \le
q+1.  \label{estimate3} \end{equation} 

Comparing the estimates
in~(\ref{estimate1}), (\ref{estimate2}), (\ref{estimate3}), we deduce that 
\begin{equation} q-16 \le |\Gamma_\pi| \le q+5,  \end{equation} 
for \emph{any} eligible pencil $\pi$. 

\subsection{} 
If $\sigma_\pi$ is the number of nonsingular members in $\pi$, then
$|\Psi_\pi| = \sigma_\pi \, (\sigma_\pi -1)$. According to 
Hirschfeld's table, $\sigma_\pi = q-2, q, q-2, q+1, q$ for classes
$(3),(14), (16), (18), (19)$ respectively. Hence $|\Psi_\pi|$ is at least
$(q-2)(q-3)$, and at most $q \, (q+1)$. We have used the former value
as the denominator of $\Upper$, and the latter value as the
denominator of $\Lower$. This gives the statement of
proposition~\ref{prop.estimate.pi}. The main theorem is now completely
proved. \qed

\subsection{} Assume that $\text{char}(q)=3$, and reconsider the
calculation in section~\ref{section.sample.calc}.  Since $f(s) =
(s+1)^2$, the quantity $\delta$  is always a perfect square. Thus 
$H_2(r,s)=0$ has a root for any $s$-value, and two roots if $s \neq
-1$. Now the same sequence of arguments shows that 
\[ \frac{|\Gamma|}{|\Psi|} \simeq \frac{2}{q} +
O\left(\frac{1}{q^2}\right). \] 
In other words, there are asymptotically twice as many (\ptc)-pairs in
characteristic $3$.

\subsection{} \label{section.ABtriangle} 
The following proposition settles an issue raised by 
definition~\ref{defn.ptc}. 
\begin{Proposition} \rm Assume that $q \ge 19$. If $(\cA, \cB)$ is a 
  (\ptc)-pair, then there exists a nondegenerate $\cA \circ \cB$ triangle. 
\end{Proposition} 
\proof 
After applying an automorphism of $\P^2$, we can assume that $\cB$ is the   
Veronese conic 
\[ \{ [1,t,t^2]: t \in \FFq \} \cup \{[0,0,1]\} \]   
defined by the equation $x \, z - y^2=0$. 
The polar line of $P = [\alpha, \beta, \gamma] \in \cA$ with  
respect to $\cB$ is $\gamma \, x - 2 \, \beta \, y + \alpha\, z=0$.  
It will intersect $\cB$ if the polynomial $\alpha \, t^2 - 2 \, \beta  
\, t + \gamma=0$ has a root in $\FFq$, i.e., if $\beta^2 - \alpha \,
\gamma$ is a square. Now choose a parametrisation 
$u \rightarrow [f_0(u), f_1(u),f_2(u)]$ of $\cA$, 
where the $f_i(u)$ are polynomials of degree at most $2$. Thus we are 
looking for solutions of the equation 
\begin{equation} v^2 - \underbrace{[f_1(u)^2 - f_0(u) \,
    f_2(u)]}_{g(u)} =0,  
\label{hasse.eqn} \end{equation} 
where $g(u)$ is a polynomial of degree at most $4$. We 
need a lower estimate on the number of solutions of this 
equation. If $g(u)$ is of degree $3$ or $4$ without repeated roots,
then~(\ref{hasse.eqn}) is an affine elliptic curve, and then it has at 
least $q - 1 -2 \, \sqrt{q}$ points by a theorem of 
Artin and Hasse (see~\cite[Ch.~V]{Silverman}). If it has repeated roots or if $\deg \,
g(u) =2$, then it is a rational curve with at least $q-2$
points as long as~(\ref{hasse.eqn}) remains irreducible. If reducible, 
then it factors into $v \pm \sqrt{g(u)}=0$, and hence must
have at least $2 \, q$ solutions. 
Since a single $u$-value leads to at most two solutions, 
in any event we have at least\footnote{Although a more refined estimate
 is possible, the increase in complexity is not worth the effort.} 
$\frac{1}{2}(q-1-2 \, \sqrt{q})$ points $P$ on
$\cA$ from which a tangent can be drawn to $\cB$. 

Now a degenerate triangle involves a common tangent to $\cA$ and 
$\cB$, of which there are at most $4$. Hence we will 
have at least one nondegenerate triangle if 
$\frac{1}{2}(q -1 - 2 \, \sqrt{q}) > 4$, 
which is assured if $q \ge 19$. \qed 

\section{A Conjecture} \label{section.higher_n}
Since Poncelet's porism is true for $n$-gons in place of triangles, it
is natural to ask whether the main theorem would 
generalise accordingly. Cayley's criterion for an arbitrary 
$n$ involves a Hankel determinant with entries taken from
the sequence $H_2, H_3, \dots$ (see~\cite{GH}). 
For instance, there exists a tetragon inscribed in
$C_r$ and circumscribed around $C_s$, if and only if 
\[ H_3 = s^6-(2r+2) \, s^5+5 \, r \, s^4-5 \, r^2 \, s^2+(2 \, r^3+2 \, r^2)
\, s-r^3=0. \]
The analogous conditions for pentagons and hexagons are respectively 
\[ \left| \begin{array}{cc} H_2 & H_3 \\ H_3 & H_4 \end{array} \right|
=0, \quad \text{and} \quad 
\left| \begin{array}{cc} H_3 & H_4 \\ H_4 & H_5 \end{array} \right|
=0, \] 
but the corresponding polynomials are already too cumbersome to write down. 
I have made some computational experiments in {\sc Maple} to count the number of
root-pairs of such polynomials; they seem to support the following conjecture: 

\begin{Conjecture} \rm 
The proportion of conic pairs in $\P^2(\FFq)$ 
satisfying the Poncelet $n$-gon condition is asymptotically equal to 
$\tau_n/q$, for some \emph{integer} value $\tau_n$. 
\end{Conjecture} 
We have $\tau_3=1$, by the main theorem of this paper. Based upon 
experimental data, the next few values are conjectured to be: 
\[ \tau_4 = 3, \quad \tau_5 = 1, \quad \tau_6 = 4, \quad 
\tau_7 = 1, \quad \tau_8 = 6, \quad \tau_9 = 2. \]

\medskip 

\centerline{--} 

\vspace{1cm}

\parbox{7cm}{ \small 
Jaydeep Chipalkatti \\
Department of Mathematics \\ 
University of Manitoba \\ 
Winnipeg, MB R3T 2N2 \\ 
Canada. \\ 
{\tt jaydeep.chipalkatti@umanitoba.ca}}
\end{document}